\documentclass[12pt]{article}
\usepackage{amsmath}
\usepackage{graphicx}
\usepackage{url}
\usepackage[
pdftex,
bookmarks=true,
colorlinks=true,
linkcolor=black,
citecolor=blue,
filecolor=black,
urlcolor=blue,
plainpages=false,
pdfpagelabels,
pdftitle={},
pdfauthor={},
pdfsubject={},
pdfkeywords={},
]{hyperref}
\usepackage[numbers,sort&compress]{natbib}

\usepackage{booktabs} 

\newcommand{\dsstrut}{\rule[-24pt]{0pt}{48pt}}

\newcommand{\mediumstrut}{\rule[-10pt]{0pt}{30pt}}

\author{J.~W.~Peterson}
\date{\today}
\title{Analytical Formulae for Two of A.~H.~Stroud's Quadrature Rules}

\newcommand{\bv}[1]{{\boldsymbol{#1}}}

\begin{document}

\maketitle

\begin{abstract}
Analytical formulae for the points and weights of two fifth-order
quadrature rules for $C_3$, the 3-cube, are given.  The rules,
originally formulated by A.~H.~Stroud in 1967, are discussed in
greater detail in terms of both the setup of the basic equations and
the method of obtaining their solutions analytically.  The primary
purpose of this paper is to better document what we feel is a
particularly practical quadrature rule (e.g.\ in finite element
calculations) and one for which we felt comprehensive information was
scarce.
\vspace{12pt}

\end{abstract}

\section{Introduction}
In 1967, A.~H.~Stroud published an article~\cite{Stroud_1967} on
fifth-degree integration formulas for several symmetric,
$n$-dimensional regions.  In the first sentence of Section 2 of that
article, Stroud mentions that ``Unless stated otherwise we assume that
$n\geq 4$.''  He goes on to give a general description of what is now
a well-known method for determining non-product quadrature rules for
several standard regions including the $n$-cube, $n$-sphere, and the
entire $n$-space.  

This description is followed by a number of tabulated quadrature rules
for specific $n$ having about six digits of precision.  A particularly
interesting rule (herein referred to as ``Stroud's first rule'') is
given for the 3-cube, ``$C_3$''.  In the last sentence of the paper,
Stroud states: ``Previously no such 13 point formula was known for
$C_3$.''  A second rule (herein referred to as ``Stroud's second
rule'') for the 3-cube, having some points outside the region, is
given as well.

Unfortunately, the generic equations for the $n$-cube given by Stroud
in~\cite{Stroud_1967} do not apply to the $n=3$ case.  This fact makes
it impossible to reproduce Stroud's tabulated results (in order to
e.g.\ compute the points and weights to a higher precision) without
investing some amount of time redoing the algebra oneself for the
$n=3$ case.  In Stroud's famous 1971 compendium on quadrature
rules~\cite{Stroud_1971}, the first rule for the 3-cube is given in
additional detail: fifteen decimal digits of precision (suitable for
double-precision calculations) are given for the points and weights,
and an eighth-order polynomial is given for determining two of the rule
parameters, but no additional guidance is given in this short table
entry.  Stroud's second rule for the 3-cube is also not given in
detail in~\cite{Stroud_1971}.

In this paper we provide additional details on Stroud's first and
second rules for the 3-cube.  In particular, we give analytical
formulae for the points and weights, which, to our knowledge, have not
been previously published.
Although we came across a number of potential references for these
rules~\cite{Cools_1993,Cools_1999,Cools_2003,Mol82,god94,Bec87}, none
were found that contained exactly the information given here.
The purpose of this paper is
therefore to provide an accessible reference for what we feel is an
especially practical rule.  For example, the analytical
formulae presented here are now being used in the general purpose
finite element library, \texttt{LibMesh}~\cite{LibMeshPaper}.
The reasons for preferring an analytical solution over
tabulated numerical values are obvious, chief among them being
adequate precision in any computing environment, even those not yet in
existence.

\section{Basic Equations for Stroud's Fifth-Order Rules for the 3-cube}
Stroud's technique for obtaining fifth-order quadrature rules for
symmetric regions $R$ is summarized as follows: choose $N$ points
$\nu_i$ and weights $A_i$ such that the approximation
\begin{equation}
  \label{eqn:rule_definition}
  \sum_{i=1}^N A_i f(\nu_i) \approx \int_R f dx
\end{equation}
is exact for all monomial functions \mbox{$f:=x_1^{\alpha_1} x_2^{\alpha_2}\!\!\ldots x_n^{\alpha_n}$} 
for which
\begin{equation}
  \label{eqn:alpha_sum}
 |\alpha| := \alpha_1 + \alpha_2 + \ldots + \alpha_n \leq 5
\end{equation}
Due to the symmetry of the regions studied, Stroud
additionally restricts the set of possible rules to only those which
contain symmetric pairs of points $\pm\nu_i$ and their negatives, both having
weight $A_i$.  (If $\nu_i$ is the origin, then its negative is not
included in the rule.)  The rules are also restricted to have
\mbox{$M=\frac{1}{2}(n^2 + n) + 1$} distinct points.  Such rules are then automatically
exact for all monomials in which $|\alpha|$ is odd.

The condition that Eqn.~\eqref{eqn:rule_definition}
be exact for monomials with \mbox{$|\alpha| = 0, 2, 4$}
is then equivalent to the matrix equation
\begin{equation}
  \label{eqn:Stroud5}
  \bv{X}^T \bv{A} \bv{X} = \bv{C}
\end{equation}
which is Stroud's Eqn.~(5), where
$
\bv{A} := \text{diag} \left(A_1, A_2, \ldots, A_M\right)
$,
$\bv{X}$ contains various (quadratic) products of the $\nu_i$, and
$\bv{C}$ contains corresponding exact monomial integrals for the given region.  (The
reader should refer to~\cite{Stroud_1967} for additional details.)  
Assuming $\bv{C}$ is non-singular, we can take inverses and thus
Eqn.~\eqref{eqn:Stroud5} is equivalent to
\begin{equation}
\label{eqn:Stroud7}
\bv{X}\bv{C}^{-1}\bv{X}^T = 
\bv{A}^{-1}
\end{equation}
which is Stroud's Eqn.~(7).  Finally, Stroud chooses the $M$ points and
weights in a special way by taking (for the particular case $n=3$, $M=7$)
\begin{equation}
  \nonumber
\left.\begin{array}{c}
\nu_1 = (\eta, \eta, \eta)
\end{array}
\right\} \text{weight } A
\end{equation}
\begin{equation}
  \nonumber
\left.\begin{array}{c}
  \nu_2 = (\lambda, \xi, \xi)  \\
  \nu_3 = (\xi, \lambda, \xi)  \\
  \nu_4 = (\xi, \xi, \lambda)  \\
\end{array}
\right\} \text{weight } B
\qquad
\left.\begin{array}{c}
  \nu_5 = (\mu, \mu, \gamma) \\
  \nu_6 = (\mu, \gamma, \mu) \\
  \nu_7 = (\gamma, \mu, \mu) \\
\end{array}\right\} \text{weight } C
\end{equation}
With this special choice of points, and taking $\eta:=0$,
Eqn.~\eqref{eqn:Stroud7} has the left-hand side
\begin{equation}
\label{eqn:xcxt}
\bv{X}\bv{C}^{-1}\bv{X}^T =
\left[
\begin{array}{c|ccc|ccc}
19   &    &\bv{M}^T_1(\lambda,\xi)  &    &     & \bv{M}^T_1(\gamma,\mu)    & \\ \hline
                                      &    &                       &    &     &            & \\
          \bv{M}_1(\lambda,\xi)     &    &\bv{M}_2(\lambda,\xi)  &    &     & \bv{M}_3   & \\
                                      &    &                       &    &     &            & \\ \hline
                                      &    &          &            &     &                        & \\
          \bv{M}_1(\gamma,\mu)      &    &\bv{M}_3  &            &     &\bv{M}_2(\gamma,\mu)    & \\
                                      &    &          &            &     &                        & 
\end{array}
\right]
\end{equation}
and the right-hand side
\begin{equation}
  \label{eqn:Amatrix}
  \bv{A}^{-1} = 32\,\text{diag}\left(A^{-1}, B^{-1}, B^{-1}, B^{-1}, C^{-1}, C^{-1}, C^{-1}\right)
\end{equation}
We have added extra horizontal and vertical lines to the matrix given in Eqn.~\eqref{eqn:xcxt} to 
emphasize the symmetric sub-structure present in the governing equations.
In Eqns.~\eqref{eqn:xcxt} and~\eqref{eqn:Amatrix} we have also introduced the following submatrices: 
\begin{equation}
  \label{eqn:m1}
  \bv{M}_1(x,y) := 
  \begin{bmatrix}
    m_1 \\
    m_1 \\ 
    m_1
  \end{bmatrix} , \qquad
  m_1(x,y) := 19 - 15x^2 - 30y^2
\end{equation}
\begin{equation}
  \nonumber 
  \bv{M}_2(x,y) := 
  \begin{bmatrix}
    m_2 & m_3 & m_3 \\ 
    m_3 & m_2 & m_3 \\
    m_3 & m_3 & m_2 \\
  \end{bmatrix},
\end{equation}
\begin{align}
  \label{eqn:m2} 
  m_2(x,y) & := 45x^{4} 
     - 30x^2 
     - 60y^2 
     + 126y^{4} 
     + 72x^2y^2 
     + 19 \\
  \label{eqn:m3} 
  m_3(x,y) & := 
     45y^{4} 
     - 30x^2 
     - 60y^2 
     + 126x^2y^2 
     + 72xy^{3} 
     + 19
\end{align}
\begin{equation}
  \nonumber 
  \bv{M}_3 := 
  \begin{bmatrix}
    m_4 & m_4 & m_5 \\ 
    m_4 & m_5 & m_4 \\
    m_5 & m_4 & m_4 \\
  \end{bmatrix},
\end{equation}
\begin{align} 
  \nonumber
  m_4  :=&
  \,45\xi^2\mu^2
     -30\xi^2
     +45\mu^2\lambda^2
     -15\lambda^2
     +45\gamma^2\xi^2 +36\lambda\xi\mu^2
     \\ \label{eqn:m4} 
     & \, 
     +36\xi^2\mu\gamma
     +36\lambda\xi\mu\gamma
     -30\mu^2-15\gamma^2+19
  \\
  \nonumber
  m_5  := &
  \,126\xi^2\mu^2
     -30\xi^2
     +45\gamma^2\lambda^2
     -15\lambda^2
     +72\lambda\xi\mu\gamma
     \\ \label{eqn:m5} 
     &\, 
     -30\mu^2-15\gamma^2+19
\end{align}

The special form assumed for the quadrature points and weights has reduced the size of
the original system of equations quite drastically, from $7^2=49$ equations to fewer than 10.

\section{Analytical Solution of the Equations}
The equations defined by $\bv{M}_1$ and $\bv{M}_2$ may be used to
determine all possible solutions for the triplets $(\lambda, \xi, B)$
and $(\gamma, \mu, C)$ independently and simultaneously.  
The following procedure is used: we rearrange Eqn.~\eqref{eqn:m1} to
solve for $y$ (resp.\ $x$) and insert it into the $m_3$
equation,~\eqref{eqn:m3}.  Since the equation for $m_3$ has odd powers
of both $x$ and $y$, we obtain two possible forms of the $m_3$
equation: one for the positive square root of $x$ (resp.\ $y$) and one
for the negative square root.  Setting the product of the positive-
and negative-root versions of Eqn.~\eqref{eqn:m3} equal to zero leads
to an eighth-order polynomial equation in $x$ (resp.~$y$) which is
equivalent to a quartic polynomial equation in $x^2$ (resp.~$y^2$).
We can solve this quartic equation for $x^2$ (resp.~$y^2$)
analytically, and, finally, the resulting $(x^2,y^2)$ pairs may be
substituted into the $m_2$ equation,~\eqref{eqn:m2} to solve for the
weights.

The eighth-order polynomial in $x$ (quartic in $x^2$) arising from the
previously-described procedure is given by
\begin{equation}
  \label{eqn:x_8th}
1330425x^{8} 
-3108780x^{6}
+2339622x^{4}
-\frac{2828796}{5}x^{2}
+361
= 0
\end{equation}
while the eighth-order equation for $y$ is
\begin{equation}
  \label{eqn:y_8th}
53217y^8
- \frac{363204}{5}y^6
+ \frac{833454}{25}y^4
- \frac{30324}{5}y^2
+ 361
= 0
\end{equation}

We can use any suitable CAS to solve Eqn.~\eqref{eqn:x_8th} and obtain the four solutions
of $x^2$ (recall that our generic $x$ variable corresponds to either $\lambda$ or $\gamma$) 
analytically.  The result is 
\begin{align}
\nonumber
x^2 & = \{ x^2_{1,2}, x^2_{3,4}\} \\
\label{eqn:x2_analytical}
    & =  \Bigg\{ \frac{1919 + 148 \sqrt {19} \pm 4 t_-}{3285}, 
                 \frac{1919 - 148 \sqrt {19} \pm 4 t_+}{3285} \Bigg\} 
\end{align}
where the short-hand notation 
\begin{equation}
    \label{eqn:t}
    t_{\pm} := \sqrt {71440 \pm 6802\sqrt {19}}
\end{equation}
is used, and will also be used throughout this paper.
In Eqn.~\eqref{eqn:x2_analytical} and those which follow, for any
quantity with an $i,j$ subscript and an ambiguous ``$\pm$'' sign, the
$i$ subscript always refers to the top sign throughout the equation,
while the $j$ subscript refers to the bottom sign.  

In a similar manner, we obtain the following analytical solutions for $y^2$ 
(recall that our generic $y$ variable corresponds to either $\xi$ or $\mu$)
by solving Eqn.~\eqref{eqn:y_8th}
\begin{align}
\nonumber
y^2 &=   \{ y^2_{1,2}, y^2_{3,4}\} \\
\label{eqn:y2_analytical}
    &= \Bigg\{ \frac{1121 - 74\sqrt{19} \mp 2 t_-}{3285},         
         \frac{1121 + 74\sqrt{19} \mp 2 t_+}{3285} \Bigg\}     
\end{align}
(Note: we have ordered the $x_i^2$ and $y_i^2$ solutions such that pairs $(x_i, y_i)$
satisfy the $m_1$ equation~\eqref{eqn:m1}.)  We can now substitute the $(x_i, y_i)$
pairs given above into the generic $m_2$ equation (repeated here)
\begin{equation}
  \label{eqn:m2_wts}
  45x^{4} 
     - 30x^2 
     - 60y^2 
     + 126y^4 
     + 72x^2y^2 
     + 19 = 32w^{-1}
\end{equation}
to solve for the generic weights $w$.  The weights $w_i$ 
(which we have already scaled by $\frac{1}{2}$, since each weight
solved for in Eqn.~\eqref{eqn:m2_wts} is actually twice the true
value) correspond to the $B$ and $C$ parameters in the
 original equations, and are given by
\begin{align}
\nonumber
w =   \{ &w_{1,2}, w_{3,4}\} \\
\nonumber
  =  \Bigg\{ & \frac{133225}{260072 + 1520\sqrt {19} \pm \left(133 + 37\sqrt {19}\right) t_-},         \\ %
\label{eqn:w_analytical}
             & \frac{133225}{260072 - 1520\sqrt {19} \pm \left(133 - 37\sqrt {19}\right) t_+} \Bigg\}    %
\end{align}
Upon simplification, the $m_5$ equation,~\eqref{eqn:m5} yields 
\begin{equation}
  \label{eqn:m5_simple}
  126\xi^2\mu^2
  +45\gamma^2\lambda^2
  +72\lambda\xi\mu\gamma
  + 19 = 0  
\end{equation}
This implies that exactly one of the four parameters $\lambda$, $\xi$,
$\mu$, $\gamma$ must be negative.  (Note: Conversely, three of the
four parameters could instead be negative, but since the rule always
includes a point and its negative, this is equivalent to one of the
four parameters being negative.)  This means that, when taking square
roots of the $x^2$ and $y^2$ values obtained previously, exactly one
negative root must be selected.  In both of his rules, Stroud~\cite{Stroud_1967}
has selected $\xi$ as the negative root, and we shall follow the same
convention here.  Using the rest of Stroud's tabulated results as a guide, we 
have compiled Table~\ref{tab:all_analytical}, which gives the corresponding
analytically-obtained $x_i$, $y_i$, and $w_i$ values.  It may of course
be verified that these analytical values also satisfy Eqns.~\eqref{eqn:m4} and \eqref{eqn:m5}.

\begin{table}[hbt]
  \begin{center}
    \caption{The middle column gives numerical approximations to the analytical
      results obtained in this section, for reference.  In the right-hand column,
      we give the corresponding parameter originally obtained by 
      Stroud~\cite{Stroud_1967}.  The 1 and 2 subscripts in the third column correspond
      to the first (with all points inside the region) and second (with some points
      outside the region) rules reported by Stroud.
      \label{tab:all_analytical}}
    \begin{tabular}[t]{lll} 
      \toprule
      $x_1$ & 1.0146309695  & $\phantom{-}\gamma_2$ \\ 
      $x_2$ & 0.7291297984  & $\phantom{-}\lambda_2$ \\ 
      $x_3$ & 0.8803044067  & $\phantom{-}\lambda_1$ \\ 
      $x_4$ & 0.0252937117  & $\phantom{-}\gamma_1$ \\ 
      \midrule
      $y_1$ & 0.3443767286  & $\phantom{-}\mu_2$ \\ 
      $y_2$ & 0.6062327951  & $-\xi_2$ \\ 
      $y_3$ & 0.4958481715  & $-\xi_1$ \\ 
      $y_4$ & 0.7956214222  & $\phantom{-}\mu_1$ \\ 
      \midrule
      $w_1$ & 0.4075948702  & $\phantom{-}C_2$ \\ 
      $w_2$ & 0.6450367090  & $\phantom{-}B_2$ \\ 
      $w_3$ & 0.5449873514  & $\phantom{-}B_1$ \\ 
      $w_4$ & 0.5076442277  & $\phantom{-}C_1$ \\ 
      \bottomrule
    \end{tabular}
  \end{center}
\end{table}

\section{Summary\label{sec:summary}}
To summarize the analytical results reported, we give both numerical
approximations and analytical forms for the various parameters in this
section.  In Tables~\ref{tab:stroudrule1}
and~\ref{tab:stroudrule1_analytical} the numerical and analytical
values, respectively, for Stroud's first quadrature rule are given.
The numerical approximations are given to 32 digits of accuracy for
convenience and because in the case of the second rule, such 
highly-accurate values have not been previously published.  
Stroud's second fifth-order rule (which is less useful for finite
element calculations due to the fact that some of the points lie
outside the region of integration) is likewise summarized in 
Tables~\ref{tab:stroudrule2} and~\ref{tab:stroudrule2_analytical}.

\begin{table}[hbt]
  \begin{center}
    \caption{Values for Stroud's first fifth order quadrature rule
      for the 3-cube, originally reported in~\cite{Stroud_1967}, 
      to 32 decimal digits.\label{tab:stroudrule1}}
    \vspace{3pt}
    \begin{tabular}{cl} 
      \toprule
      $\eta$    & $\phantom{-}0.00000000000000000000000000000000\text{E}\!+\!00$ \\ 
      $\lambda$ & $\phantom{-}8.80304406699309780477378182098603\text{E}\!-\!01$  \\
      $\xi$     & $-4.95848171425711152814212423642879\text{E}\!-\!01$ \\
      $\mu$     & $\phantom{-}7.95621422164095415429824825675787\text{E}\!-\!01$ \\
      $\gamma$  & $\phantom{-}2.52937117448425813473892559293236\text{E}\!-\!02$ \\
      \midrule
      $A$ & $\phantom{-}1.68421052631578947368421052631579\text{E}\!+\!00$ \\
      $B$ & $\phantom{-}5.44987351277576716846907821808944\text{E}\!-\!01$ \\ 
      $C$ & $\phantom{-}5.07644227669791704205723757138424\text{E}\!-\!01$ \\ 
      \bottomrule
    \end{tabular}
  \end{center}
\end{table}

%
\begin{table}[hbt]
  \begin{center}
    \caption{Analytical representations for Stroud's first fifth order
      quadrature rule for the 3-cube, originally reported
      in~\cite{Stroud_1967}.  The results correspond to the numerical
      approximations given in Table~\ref{tab:stroudrule1}.  When taking the square root,
      the positive root is always assumed unless specified otherwise.  
      \label{tab:stroudrule1_analytical}}
  \end{center}
  \vspace{-3pt}
    %
    %
    \begin{tabular}[t]{|c|l|} \hline
      $\eta$    & $0$ \dsstrut \\ 
      $\lambda$ & $\phantom{-}\displaystyle \sqrt{\frac{1919 - 148 \sqrt {19} + 4 t_+}{3285}}$ \dsstrut \\
      $\xi$     & $-\displaystyle \sqrt{\frac{1121 + 74\sqrt{19} - 2 t_+}{3285}}$ \dsstrut \\
      $\mu$     & $\phantom{-}\displaystyle \sqrt{\frac{1121 + 74\sqrt{19} + 2 t_+}{3285}}$ \dsstrut \\ 
      $\gamma$  & $\phantom{-}\displaystyle \sqrt{\frac{1919 - 148 \sqrt {19} - 4 t_+}{3285}}$ \dsstrut \\ \hline
    \end{tabular}
    \begin{tabular}[t]{|c|l|} \hline
      $A$ & $\displaystyle \frac{32}{19}$ \dsstrut \\
      $B$ & $\displaystyle \frac{133225}{260072 - 1520\sqrt{19} + \left(133 - 37\sqrt{19}\right) t_+}$ \dsstrut \\ 
      $C$ & $\displaystyle \frac{133225}{260072 - 1520\sqrt{19} 
                           - \left(133 - 37\sqrt{19}\right) t_+}$ \dsstrut \\ \hline\hline
      \multicolumn{2}{|c|}{} \\ 
      \multicolumn{2}{|c|}{$\displaystyle t_{\pm} := \displaystyle \sqrt {71440 \pm 6802\sqrt {19}}$ \mediumstrut} \\                   
      \multicolumn{2}{|c|}{} \\ 
      \hline
    \end{tabular}
\end{table}

\begin{table}[hbt]
  \begin{center}
    \caption{More accurate values for Stroud's second fifth order quadrature rule
      for the 3-cube, originally reported in~\cite{Stroud_1967}, to 32 decimal 
      digits.\label{tab:stroudrule2}}
    \vspace{3pt}
    \begin{tabular}{cl} 
      \toprule
      $\eta$    & $\phantom{-}0.00000000000000000000000000000000\text{E}\!+\!00$ \\ 
      $\lambda$ & $\phantom{-}7.29129798350178619517304350180274\text{E}\!-\!01$ \\
      $\xi$     & $          -6.06232795147414690867819787871004\text{E}\!-\!01$ \\
      $\mu$     & $\phantom{-}3.44376728634554308789703229927940\text{E}\!-\!01$  \\
      $\gamma$  & $\phantom{-}1.01463096947441524044348381753366\text{E}\!+\!00$ \\
      \midrule
      $A$ & $\phantom{-}1.68421052631578947368421052631579\text{E}\!+\!00$ \\
      $B$ & $\phantom{-}6.45036708927015064146303620658267\text{E}\!-\!01$ \\ 
      $C$ & $\phantom{-}4.07594870020353356906327958289101\text{E}\!-\!01$ \\ 
      \bottomrule
    \end{tabular}
  \end{center}
\end{table}

\begin{table}[hbt]
  \begin{center}
    \caption{Analytical representations for Stroud's second fifth-order
      quadrature rule for the 3-cube, originally reported
      in~\cite{Stroud_1967}.  The results correspond to the numerical
      approximations given in Table~\ref{tab:stroudrule2}.  When taking the square root,
      the positive root is always assumed unless specified otherwise.  
      \label{tab:stroudrule2_analytical}}
  \end{center}
  \vspace{-3pt}
    %
    %
  \begin{tabular}[t]{|c|l|} \hline
    $\eta$    & $0$ \dsstrut \\ 
    $\lambda$ & $\phantom{-}\displaystyle\sqrt{\frac{1919 + 148 \sqrt {19} - 4 t_-}{3285}}$ \dsstrut \\
    $\xi$     & $-\displaystyle\sqrt{\frac{1121 - 74\sqrt{19} + 2 t_-}{3285}}$ \dsstrut \\
    $\mu$     & $\phantom{-}\displaystyle \sqrt{\frac{1121 - 74\sqrt{19} - 2 t_{-}}{3285}}$ \dsstrut \\ 
    $\gamma$  & $\phantom{-}\displaystyle \sqrt{\frac{1919 + 148 \sqrt {19} + 4 t_{-}}{3285}}$ \dsstrut \\ \hline 
  \end{tabular}
  \begin{tabular}[t]{|c|l|} \hline
    $A$ & $\displaystyle \frac{32}{19}$ \dsstrut \\
    $B$ & $\displaystyle 
    \frac{133225}{260072 + 1520\sqrt{19} - \left(133 + 37\sqrt{19}\right) t_-}$ \dsstrut \\ 
    $C$ & $\displaystyle \frac{133225}{260072 + 1520\sqrt{19} 
      + \left(133 + 37\sqrt{19}\right) t_{-}}$ \dsstrut \\ \hline\hline
      \multicolumn{2}{|c|}{} \\ 
      \multicolumn{2}{|c|}{$\displaystyle t_{\pm} := \displaystyle \sqrt {71440 \pm 6802\sqrt {19}}$ \mediumstrut} \\                   
      \multicolumn{2}{|c|}{} \\ 
      \hline
  \end{tabular}
\end{table}

\clearpage 
\bibliographystyle{plain}  
\bibliography{stroud_hex13}

\end{document}